\documentclass[number,citesort,MSNbibl,seceqn,dvips]{arxbj}
\usepackage{subenv,upgreek}
\usepackage{mathrsfs}

\aid{0}
\volume{18}
\issue{1}
\pubyear{2012}
\firstpage{46}
\lastpage{63}
\doi{10.3150/10-BEJ329}

\makeatletter
\newcommand{\coloneqq}{:=}
\newcommand{\eqqcolon}{=:}
\newcommand{\eqref}[1]{(\ref{#1})}
\newcommand{\text}{\mathrm}
\newcommand{\cancel}{\mathop{\,\mbox{\parbox[c][5.5pt][b]{11.5pt}{\fontsize{11.5}{11.5}\selectfont{$\not$}}}\!\!\!\!\!\!\!\!\!}}

\newcommand{\Cov}{\operatorname{Cov}}
\newcommand{\Var}{\operatorname{Var}}
\newcommand{\rank}{\operatorname{rank}}
\newcommand{\adj}{\operatorname{adj}}
\newcommand{\re}{\operatorname{Re}}
\newcommand{\imag}{\operatorname{Im}}
\newcommand{\E}{\mathbb{E}}
\newcommand{\Pb}{\mathbb{P}}
\newcommand{\I}{\mathbb{I}}
\newcommand{\N}{\mathbb{N}}
\newcommand{\Z}{\mathbb{Z}}
\newcommand{\R}{\mathbb{R}}
\newcommand{\C}{\mathbb{C}}
\newcommand{\Y}{\mathbf{Y}}
\newcommand{\X}{\mathbf{X}}
\newcommand{\Lb}{\mathbf{L}}
\newcommand{\J}{\mathbf{J}}
\newcommand{\W}{\mathbf{W}}
\newcommand{\V}{\mathbf{V}}
\newcommand{\ZZ}{\mathbf{Z}}
\newcommand{\NN}{\mathbf{N}}
\newcommand{\M}{\mathbf{M}}
\newcommand{\beps}{\bolds{\varepsilon}}
\newcommand{\bepsh}{\bolds{\varepsilon}^{(h)}}
\newcommand{\bu}{\mathbf{u}}
\newcommand{\bx}{\mathbf{x}}
\newcommand{\bxi}{\bolds{\xi}}
\newcommand{\bb}{\mathbf{b}}
\newcommand{\bs}{\mathbf{s}}
\newcommand{\bzero}{\mathbf{0}}
\newcommand{\be}{\mathbf{e}}
\newcommand{\vecM}{\mathscr{M}}
\newcommand{\nuL}{\nu^{\Lb}}

\newtheorem{theorem}{Theorem}[section]
\newproclaim{definition}[theorem]{Definition}
\newtheorem{proposition}[theorem]{Proposition}
\newtheorem{corollary}[theorem]{Corollary}
\newtheorem{lemma}[theorem]{Lemma}
\newproclaim{assumptionlevy}{Assumption}

\newproclaim{assumptioneigen}{Assumption}

\makeatother

\begin{document}
\begin{frontmatter}

\title{Multivariate CARMA processes, continuous-time state space models and complete regularity of the innovations of the sampled processes}
\runtitle{Complete regularity of the innovations of sampled MCARMA processes}

\begin{aug}
\author{\fnms{Eckhard} \snm{Schlemm}\corref{}\ead[label=e1,mark]{schlemm@ma.tum.de}\ead[label=u1,url]{http://www-m4.ma.tum.de}} \and
\author{\fnms{Robert} \snm{Stelzer}\ead[label=e2,mark]{rstelzer@ma.tum.de}\ead[label=u1,url]{http://www-m4.ma.tum.de}}
\runauthor{E. Schlemm and R. Stelzer}
\address{TUM Institute for Advanced Study and Zentrum Mathematik,
Technische Universit\"{a}t M\"{u}nchen,\\
Boltzmannstr.~3,
85748 Garching, Germany. \printead{e1,e2};\\ \printead{u1}}
\end{aug}

\received{\smonth{5} \syear{2010}}
\revised{\smonth{9} \syear{2010}}

%
\begin{abstract}
The class of multivariate L\'{e}vy-driven autoregressive moving average
(MCARMA) processes, the continu\-ous-time analogs of the classical vector
ARMA processes, is shown to be equivalent to the class of
continuous-time state space models. The linear innovations of the weak
ARMA process arising from sampling an MCARMA process at an equidistant
grid are proved to be exponentially completely regular ($\beta$-mixing)
under a mild continuity assumption on the driving L\'{e}vy process. It
is verified that this continuity assumption is satisfied in most
practically relevant situations, including the case where the driving
L\'{e}vy process has a non-singular Gaussian component, is compound
Poisson with an absolutely continuous jump size distribution or has an
infinite L\'{e}vy measure admitting a density around zero.
\end{abstract}

%
\begin{keyword}
\kwd{complete regularity}
\kwd{linear innovations}
\kwd{multivariate CARMA process}
\kwd{sampling}
\kwd{state space representation}
\kwd{strong mixing}
\kwd{vector ARMA process}
\end{keyword}

\end{frontmatter}

\section{Introduction}
CARMA processes are the continuous-time analogs of the widely known
discrete-time ARMA processes (see, e.g., \cite{brockwell1991tst} for
a comprehensive introduction); they were first defined in \cite
{doob1944egp} in the univariate Gaussian setting and have stimulated a
considerable amount of research in recent years (see, e.g., \cite
{brockwell2001cta} and references therein). In particular, the
restriction of the driving process to Brownian motion was relaxed and
\cite{brockwell2001ldc} allowed for L\'{e}vy processes with finite
logarithmic moments. Because of their applicability to irregularly
spaced observations and high-frequency data, they have turned out to be
a versatile and powerful tool in the modeling of phenomena from the
natural sciences, engineering and finance. Recently, \cite
{marquardt2007multivariate} extended the concept to multivariate CARMA
(MCARMA) processes with the intention of being able to model the joint
behavior of several dependent time series.
MCARMA processes are thus the continuous-time analogs of discrete-time
vector ARMA (VARMA) models (see, e.g., \cite{lutkepohl2005nim}).

The aim of this paper is twofold: first, we establish the equivalence
between MCARMA and multivariate continuous-time state space models, a
correspondence which is well known in the discrete-time setting \cite
{hannan1987stl}; second, we investigate the probabilistic properties
of the discrete-time process obtained by recording the values of an
MCARMA process at discrete, equally spaced points in time. A detailed
understanding of the innovations of the weak VARMA process which arises
is a prerequisite for proving asymptotic properties of statistics of a
discretely observed MCARMA process. One notion of asymptotic
independence which is very useful in this context is complete
regularity (see Section \ref{section-results} for a~precise definition)
and we show that the innovations of a discretized MCARMA process have
this desirable property. Our results therefore not only provide
important insight into the probabilistic structure of CARMA processes,
but they are also fundamental to the development of an estimation
theory for non-Gaussian continuous-time state space models based on
equidistant observations.

In this paper, we show that a sampled MCARMA process is a discrete-time
VARMA process with dependent innovations. While the mixing behavior of
ARMA and more general linear processes is fairly well understood (see,
e.g., \cite{athreya1986mixing,mokkadem1988mixing,pham1985some}), the
mixing properties of the innovations of a sampled continuous-time
process have received very little attention. From \cite
{brockwell2009existence}, it is only known that the innovations of a
discretized univariate L\'{e}vy-driven CARMA process are weak white
noise, which, by itself, is typically of little help in applications.
We show that the linear innovations of a sampled MCARMA process satisfy
a set of VARMA equations and we conclude that under a mild continuity
assumption on the driving L\'{e}vy process, they are geometrically
completely regular and, in particular, geometrically strongly mixing.
This continuity assumption is further shown to be satisfied for most of
the practically relevant choices of the driving L\'{e}vy process,
including processes with a non-singular Gaussian component, as well as
compound Poisson processes with an absolutely continuous jump size
distribution and infinite activity processes whose L\'{e}vy measures
admit a density in a neighborhood of zero.

This paper is structured as follows. In Section \ref{section-MLP} we
review some well-known properties of L\'{e}vy processes, which we will
use later. The class of multivariate CARMA processes, in a slightly
more general form than in the original definition of \cite
{marquardt2007multivariate}, is described in detail in Section~\ref
{section-MCARMA} and shown to be equivalent to the class of
continuous-time state space models. In Section~\ref{section-results}
the main result about the mixing properties of the sampled processes is
stated and demonstrated to be applicable in many practical situations.
The proofs of the results are presented in Section~\ref{sec5}.

We use the following notation. The space of $m\times n$ matrices with
entries in the ring $\mathbb{K}$ is denoted by $M_{m,n}(\mathbb{K})$ or
$M_{m}(\mathbb{K})$ if $m=n$. $A^{\mathrm{T}}$ denotes the transpose of the matrix
$A$, the matrices $\I_m$ and $0_m$ are the identity and the zero
element of $M_m(\mathbb{K})$, respectively, and $A \otimes B$ stands
for the Kronecker product of the matrices $A$ and $B$. The zero vector
in $\R^m$ is denoted by $\bzero_m$, and $\|\cdot\|$ and $\langle
\cdot
,\cdot\rangle$ represent the Euclidean norm and inner product,
respectively. Finally, $\mathbb{K}[z]$ ($\mathbb{K}\{z\}$) is the ring
of polynomial (rational) expressions in $z$ over $\mathbb{K}$ and
$I_B(\cdot)$ is the indicator function of the set $B$.

\section{Multivariate L\'{e}vy processes}
\vspace*{-2pt}
\label{section-MLP}
In this section we review the definition of a multivariate L\'{e}vy
process and some elementary facts about these processes which we will
use later. More details and proofs can be found in, for instance, \cite
{sato1999lpa}.\vspace*{-2pt}
\begin{definition}
A (one-sided) $\R^m$-valued \textit{L\'{e}vy process} $(\Lb
(t))_{t\geq0}$
is a stochastic process with stationary, independent increments,
continuous in probability and satisfying $\Lb(0)=\bzero_m$ almost
surely.\vspace*{-2pt}
\end{definition}

Every $\R^m$-valued L\'{e}vy process $(\Lb(t))_{t\geq0}$ can be assumed
to be c\`{a}dl\`{a}g and is completely characterized by its characteristic
function in the L\'{e}vy--Khintchine form $\E\mathrm{e}^{\mathrm
{i}\langle\bu,\Lb
(t)\rangle}=\exp\{t\psi^{\Lb}(\bu)\}$, $\bu\in\R^m$, $t\geq0$, where
$\psi^{\Lb}$ has the special form
\[
\label{characfunc}
\psi^{\Lb}(\bu)=\mathrm{i}\langle\bolds{\gamma},\bu\rangle
-\frac
{1}{2}\langle\bu,\Sigma^{\mathcal{G}}\bu\rangle+\int_{\R^m}
\bigl[\mathrm{e}^{\mathrm{i}\langle\bu,\bx\rangle}-1-\mathrm
{i}\langle\bu,\bx\rangle I_{\{\|x\|\leq1\}
} \bigr]\nuL(\mathrm{d}\bx).
\]
The vector $\gamma\in\R^m$ is called the \textit{drift}, the non-negative
definite, symmetric $m\times m$ matrix~$\Sigma^{\mathcal{G}}$ is the
\textit{Gaussian covariance matrix} and $\nuL$ is a measure on $\R^m$,
referred to as the \textit{L\'{e}vy measure}, satisfying
\[
\nuL(\{\bzero_m\})=0,\qquad\int_{\R^m}\min(\|\bx\|^2,1)\nuL
(\mathrm
{d}\bx)<\infty.
\]
We will work with two-sided L\'{e}vy processes $\Lb=(\Lb(t))_{t\in\R}$.
These are obtained from two independent copies $(\Lb_1(t))_{t\geq0}$,
$(\Lb_2(t))_{t\geq0}$ of a one-sided L\'{e}vy process via the construction
\[
\Lb(t)=
\cases{\displaystyle
\Lb_1(t), & \quad $t\geq0$, \cr\displaystyle
-\lim_{s\nearrow-t}\Lb_2(s), & \quad $t<0$.
}
\]
Throughout the paper, we restrict our attention to L\'{e}vy processes
with zero means and finite second moments.\vspace*{-2pt}
\begin{assumptionlevy}
\label{assumlevy2mom}
The L\'{e}vy process $\Lb$ satisfies $\E\Lb(1)=0$ and $\E\|\Lb
(1)\|^2<\infty$.\vspace*{-2pt}
\end{assumptionlevy}

The assumption $\E\Lb(1)=0$ is made only for notational convenience and
is not essential for our results to hold. The premise that $\Lb$ has
finite variance is, in contrast, a true restriction, which is very
often made in the analysis of (C)ARMA processes. The treatment of the
infinite variance case requires different techniques and often does not
lead to comparable results. It is well known that $\Lb$ has finite
second moments if and only if $\int_{\|\bx\|\geq1}{\|\bx\|^2\nu
(\mathrm
{d}\bx)}$ is finite, and that $\Sigma^{\Lb}=\E\Lb(1)\Lb(1)^{\mathrm{T}}$ is then
given by $\int_{\mathbb{R}^m}{\bx\bx^{\mathrm{T}}}\nuL(\mathrm{d}\bx)+\Sigma^{\mathcal{G}}$.

\vspace*{-2pt}\section{MCARMA processes and state space models}\vspace*{-2pt}
\label{section-MCARMA}
If $\Lb$ is a two-sided L\'{e}vy process with values in $\R^m$ and
$p>q$ are positive integers, then the $d$-dimensional\vadjust{\goodbreak} $\Lb$-driven
autoregressive moving average (\textit{MCARMA}) process with
autoregressive polynomial
%
%
\begin{subequation}
\label{ARMApoly}
\begin{equation}
\label{ARpoly}
z\mapsto P(z)\coloneqq\I_dz^p+A_1z^{p-1}+\cdots +A_p\in M_d(\R[z])
\end{equation}
and moving average polynomial
\begin{equation}
\label{MApoly}
z\mapsto Q(z)\coloneqq B_0z^q+B_1z^{q-1}+\cdots +B_q\in M_{d,m}(\R[z])
\end{equation}
\end{subequation}
is thought of as the solution to the formal differential equation
%
%
\begin{equation}
\label{formalDE}
P(D)\Y(t)=Q(D)D\Lb(t),\qquad D\equiv\frac{\mathrm{d}}{\mathrm{d}t},
\end{equation}
which is the continuous-time analog of the discrete-time ARMA
equations. We note that we allow for the driving L\'{e}vy process $\Lb$
and the $\Lb$-driven MCARMA process to have different dimensions and
thus slightly extend the original definition of \cite
{marquardt2007multivariate}. All the results we need from \cite
{marquardt2007multivariate} are easily seen to continue to hold in this
more general setting. Since, in general, L\'{e}vy processes are not
differentiable, equation \eqref{formalDE} is purely formal and, as
usual, interpreted as being equivalent to the \textit{state space
representation}
%
%
\begin{equation}
\label{MCARMAssm}
\mathrm{d}\mathbf{G}(t)=\mathcal{A}\mathbf{G}(t)\,\mathrm
{d}t+\mathcal{B} \,\mathrm{d}\Lb(t),\qquad\Y(t)= \mathcal
{C}\mathbf
{G}(t),\qquad t\in\R,
\end{equation}
where $\mathcal{A},\mathcal{B}, \mathcal{C}$ are given by
%
%
\begin{subequation}
\label{MCARMAcoeffABC}
\begin{eqnarray}
\label{MCARMAcoeffA} \mathcal{A} &=& \pmatrix{\displaystyle
0 & \I_d & 0 & \ldots& 0 \cr\displaystyle
0 & 0 & \I_d & \ddots& \vdots\cr\displaystyle
\vdots&& \ddots& \ddots& 0\cr\displaystyle
0 & \ldots& \ldots& 0 & \I_d\cr\displaystyle
-A_p & -A_{p-1} & \ldots& \ldots& -A_1
}
\in M_{pd}(\R),\\
\label{MCARMAcoeffB}\mathcal{B}&=& \pmatrix{\displaystyle\beta_1^{\mathrm{T}}
& \cdots& \beta_p^{\mathrm{T}}
}
^{\mathrm{T}}\in M_{pd,m}(\R),\nonumber
\\[-8pt]
\\[-8pt]
&&\hspace*{4pt}  \beta_{p-j} = -I_{\{0,\ldots,q\}}(j)
\Biggl[\sum_{i=1}^{p-j-1}{A_i\beta_{p-j-i}+B_{q-j}} \Biggr]
\nonumber
\end{eqnarray}
and
\begin{equation}
\label{MCARMAcoeffC}
\mathcal{C}= (\I_d,0_d,\ldots,0_d )\in
M_{d,pd}(\R).
\end{equation}
\end{subequation}
In view of representation \eqref{MCARMAssm}, MCARMA processes are
linear continuous-time state space models. We will consider this class
of processes and see that it is in fact equivalent to the class of
MCARMA models.
\begin{definition}
An $\R^d$-valued \textit{continuous-time linear state space model}
$(A,B,C,\Lb)$ of dimension $N$ is characterized by an $\R^m$-valued
driving L\'{e}vy process $\Lb$, a state transition matrix $A\in M_N(\R
)$, an input matrix $B\in M_{N,m}(\R)$ and an observation matrix $C\in
M_{d,N}(\R)$. It consists of a state equation of Ornstein--Uhlenbeck type
%
%
\begin{subequation}
\label{ssm}
\begin{equation}
\label{stateeq}
\mathrm{d}\X(t) =A\X(t) \,\mathrm{d}t+B \,\mathrm{d}\Lb(t)\vadjust{\goodbreak}
\end{equation}
and an observation equation
\begin{equation}
\label{obseq}
\Y(t) = C\X(t).
\end{equation}
\end{subequation}
The $\R^N$-valued process $\X=(\X(t))_{t\in\R}$ is the \textit{state
vector process} and $\Y=(\Y(t))_{t\in\R}$ is the \textit{output process}.
\end{definition}

A solution $\Y$ to equations \eqref{ssm} is called \textit{causal} if for all $t$,
$\Y(t)$ is independent of the $\sigma$-algebra generated by $\{\Lb
(s)\dvtx s>t\}$. Every solution to equation \eqref{stateeq} satisfies
%
%
\begin{equation}
\label{markovstate}
\X(t) = \mathrm{e}^{A(t-s)}\X(s)+\int_s^t{\mathrm{e}^{A(t-u)}B\,
\mathrm{d}\Lb(u)},\qquad
s,t\in\R,\  s<t.
\end{equation}
The independent increment property of L\'{e}vy processes implies that
$\X$ is a Markov process. We always work under the following standard
causal stationarity assumption.
\begin{assumptioneigen}
\label{assumEigen1}
The eigenvalues of $A$ have strictly negative real parts.
\end{assumptioneigen}

The following is well known \cite{sato1984operator} and recalls
conditions for the existence of a stationary causal solution of the
state equation \eqref{stateeq} for easy reference.
\begin{proposition}
If Assumptions \textup{\ref{assumlevy2mom}} and \textup{\ref{assumEigen1}} hold, then equation \textup{\eqref{stateeq}}
has a unique strictly stationary, causal solution $\X$ given by
%
%
\begin{equation}
\X(t)=\int_{-\infty}^t{\mathrm{e}^{A(t-u)}B\,\mathrm{d}\Lb
(u)},\qquad t\in\R,
\end{equation}
which has the same distribution as $\int_0^\infty{\mathrm{e}^{Au}B\,
\mathrm{d}\Lb
(u)}$. Moreover, $\X(t)$ has mean zero,
%
%
\begin{subequation}
\begin{eqnarray}
\label{Varssm}\Var(\X(t))&=&\E\X(t)\X(t)^{\mathrm{T}}\eqqcolon\Gamma
_0=\int
_0^\infty{\mathrm{e}^{Au}B\Sigma^{\Lb}B^{\mathrm{T}}\mathrm{e}^{A^\mathrm{T}u}}\,
\mathrm{d}u,\\
\label{Covssm}\Cov\bigl(\X(t+h),\X(t)\bigr)&=&\E\X(t+h)\X(t)^{\mathrm{T}}
= \mathrm{e}^{Ah}\Gamma
_0,\qquad h\geq0,
\end{eqnarray}
\end{subequation}
and $\Gamma_0$ satisfies $A\Gamma_0+\Gamma_0A^{\mathrm{T}} = -B\Sigma^{\Lb}B^{\mathrm{T}}$.
\end{proposition}

It is an immediate consequence that the output process $\Y$ has mean
zero and autocovariance function $h\mapsto\gamma_{\Y}(h) = C\mathrm
{e}^{Ah}\Gamma
_0C^{\mathrm{T}}$, and that $\Y$ can be written as a moving average of the driving
L\'{e}vy process as
%
%
\begin{equation}
\label{ssmMA}
\Y(t)=\int_{-\infty}^\infty{g(t-u)\,\mathrm{d}\Lb(u)},\qquad t\in
\R;\qquad
g(t)=C\mathrm{e}^{At}BI_{[0,\infty)}(t).
\end{equation}
These equations serve, with $A$, $B$ and $C$ defined as in equations
\eqref{MCARMAcoeffABC}, as the definition of an MCARMA process with
autoregressive and moving average polynomials given by equations \eqref
{ARMApoly}. It shows that the behavior of the process $\Y$ depends on
the values of the individual matrices $A,B,C$ only through the products
$C\mathrm{e}^{At}B$, $t\in\R$. These products are, in turn,
intimately related\vadjust{\goodbreak}
to the rational matrix function $H\dvtx z\mapsto C(z\I_N-A)^{-1}B$,
which is
called the \textit{transfer function} of the state space model \eqref
{ssm}. A pair $(P,Q)$, $P\in M_d(\R[z])$, $Q\in M_{d.m}(\R[z])$, of
rational matrix functions is a \textit{left matrix fraction description}
for the rational matrix function $H\in M_d(\R\{z\})$ if
$P(z)^{-1}Q(z)=H(z)$ for all $z\in\C$. The next theorem gives an answer
to the question of what other state space representations besides~\eqref
{MCARMAssm} can be used to define an MCARMA process. The proof is given
in Section~\ref{sec5}.

\begin{theorem}
\label{equivalenceMCARMASSM}
If $(P,Q)$ is a left matrix fraction description for the transfer
function $z\mapsto C(z\I_N-A)^{-1}B$, then the stationary solution $\Y$
of the state space model $(A,B,C,\Lb)$ defined by equations \eqref{ssm} is an
$\Lb
$-driven MCARMA process with autoregressive polynomial $P$ and moving
average polynomial $Q$.
\end{theorem}

\begin{corollary}
The classes of MCARMA and causal continuous-time state space models are
equivalent.
\end{corollary}

\begin{pf}
By definition, every MCARMA process is the output process of a state
space model. Conversely, given any state space model $(A,B,C,\Lb)$ with
output process $\Y$,~\cite{caines1988linear}, Appendix~2, Theorem 8,
shows that the transfer function $H\dvtx z\mapsto C(zI_N-A)^{-1}B$ possesses
a left matrix fraction description $H(z)=P(z)^{-1}Q(z)$. Hence, by
Theorem \ref{equivalenceMCARMASSM}, $\Y$~is an MCARMA process.
\noqed\mbox{}\qed
\end{pf}

\section{Complete regularity of the innovations of a~sampled MCARMA process}
\label{section-results}
For a continuous-time stochastic process $\Y=(\Y(t))_{t\in\R}$ and a
positive constant $h$, the corresponding sampled process $\Y^{(h)}=(\Y
^{(h)}_n)_{n\in\mathbb{Z}}$ is defined by $\Y^{(h)}_n=\Y(nh)$. A common
problem in applications is the estimation of a set of model parameters
based on observations of the values of a realization of a
continuous-time process at equally spaced points in time. In order to
make MCARMA processes amenable to parameter inference from
equidistantly sampled observations, it is important to have a good
understanding of the probabilistic properties of $\Y^{(h)}$. One such
property which has turned out to be useful for the derivation of
asymptotic properties of estimators is \textit{mixing}, for which there
are several different notions (see, e.g., \cite
{bradley2007introduction} for a detailed exposition). Let $I$ denote
$\Z$ or~$\R$. For a stationary process $\X=(X_n)_{n\in I}$ on some
probability space $(\Omega,\mathscr{F},\Pb)$, we write $\mathscr{F}_n^m=
\sigma(X_j\dvtx j\in I, n < j < m)$, $-\infty\leq n<m\leq\infty$. The
$\alpha
$-mixing coefficients $(\alpha(m))_{m\in I}$ are then defined by
\[
\alpha(m) = \sup_{A\in\mathscr{F}_{-\infty}^0,B\in\mathscr
{F}_m^\infty
}|\Pb(A\cap B)-\Pb(A)\Pb(B)|.
\]
If $\lim_{m\to\infty}\alpha(m)=0$, then the process $\X$ is called
\textit{strongly mixing}, and if there exist constants $C>0$ and
$0<\lambda<1$
such that $\alpha_m<C\lambda^m$, $m\geq1$, it is called \textit
{exponentially strongly mixing}. The $\beta$-mixing coefficients
$(\beta
(m))_{m\in I}$ are similarly defined as
\[
\beta(m) = \E\sup_{B\in\mathscr{F}_m^\infty}\bigl|\Pb(B|\mathscr
{F}_{-\infty
}^0)-\Pb(B)\bigr|.
\]
If $\lim_{m\to\infty}\beta(m)=0$, then the process $\X$ is called
\textit{completely regular} or \textit{$\beta$-mixing}, and if there exist
constants $C>0$ and $0<\lambda<1$ such that $\beta_m<C\lambda^m$,
$m\geq1$, it is called \textit{exponentially completely regular}. It is
clear from these definitions that $\alpha(m)\leq\beta(m)$ and that
(exponential) complete regularity implies (exponential) strong mixing.
It has been shown in \cite{marquardt2007multivariate}, Proposition~3.34, that every causal MCARMA process $\Y$
with a finite $\kappa$th moment, $\kappa>0$, is strongly mixing and
this naturally carries over to the sampled~pro\-cess~$\Y^{(h)}$. In this
paper, we therefore do not investigate the mixing properties of the
process~$\Y^{(h)}$ itself, but rather of its linear
innovations.\vspace*{-1pt}
\begin{definition}
Let $(\Y_n)_{n\in\Z}$ be an $\R^d$-valued stationary stochastic process
with finite second moments. The \textit{linear innovations} $(\beps
_n)_{n\in\Z}$ of $(\Y_n)_{n\in\Z}$ are then defined by
%
%
\begin{equation}
\label{DefInno}
\beps_n=\Y_n-P_{n-1}\Y_n,\qquad P_n=\mbox{orthogonal projection onto }
\overline{\operatorname{span}} \{\Y_\nu\dvtx-\infty<\nu\leq n \},
\end{equation}
where the closure is taken in the Hilbert space of square-integrable
random variables with inner product $(X,Y)\mapsto\E\langle X,Y\rangle$.\vspace*{-1pt}
\end{definition}

From now on, we work under an additional assumption, which is standard
in the univariate case.\vspace*{-1pt}
\begin{assumptioneigen}
\label{assumEigen2}
The eigenvalues $\lambda_1,\ldots,\lambda_N$ of the state transition
matrix $A$ in equation \textup{\eqref{stateeq}} are distinct.\vspace*{-1pt}
\end{assumptioneigen}

A polynomial $p\in M_d(\C[z])$ is called \textit{monic} if its leading
coefficient is equal to $\I_d$ and \textit{Schur-stable} if the zeros of
$z\mapsto\det p(z)$ all lie in the complement of the closed unit disc.
We first give a semi-explicit construction of a weak VARMA
representation of $\Y^{(h)}$ with complex-valued coefficient matrices,
a generalization of \cite{brockwell2010estimation}, Proposition 3.\vspace*{-1pt}
\begin{theorem}
\label{weakARMAthm}
Assume that $\Y$ is the output process of the state space system
\textup{\eqref{ssm}} satisfying Assumptions  \textup{\ref{assumlevy2mom}}, \textup{\ref{assumEigen1}}, \textup{\ref{assumEigen2}}, and $\Y^{(h)}$ is its sampled version with linear
innovations~$\bepsh$. Define the Schur-stable polynomial $\varphi\in
\C
[z]$ by
%
%
\begin{equation}
\varphi(z)=\prod_{\nu=1}^N (1-\mathrm{e}^{h\lambda_\nu}z
)\eqqcolon
(1-\varphi_1 z-\cdots-\varphi_N z^N ).
\end{equation}
There then exists a monic Schur-stable polynomial $\Theta\in M_d(\C
[z])$ of degree at most $N-1$ such that
%
%
\begin{equation}
\label{weakARMA}
\varphi(B)\Y_n^{(h)}=\Theta(B)\bepsh_n,\qquad n\in\Z,
\end{equation}
where $B$ denotes the backshift operator, that is, $B^j\Y^{(h)}_n=\Y
^{(h)}_{n-j}$ for every non-negative integer~$j$.\vspace*{-1pt}
\end{theorem}

This result is very important for the proof of the mixing properties of
the innovations sequence $\bepsh$ because it establishes an explicit
linear relationship between $\bepsh$ and $\Y^{(h)}$. A good
understanding of the mixing properties of $\bepsh$ is not only
theoretically interesting, but is also\vadjust{\goodbreak} practically of considerable
relevance for the purpose of statistical inference for multivariate
CARMA processes. One estimation procedure in which the importance of
the mixing properties of the innovations of the sampled process is
clearly visible is Gaussian maximum likelihoood (GML) estimation.
Assume that $\Theta\subset\R^s$ is a compact parameter set and that a
parametric family of MCARMA processes is given by the mapping $\Theta
\ni
\vartheta\mapsto(A_{\vartheta},B_{\vartheta},C_{\vartheta},\Lb
_{\vartheta})$. It follows from Theorem \ref{weakARMAthm} and \cite{brockwell1991tst}, Section 11.5, that the Gaussian likelihood of
observations $\overline{\mathbf{y}}^L=(\mathbf{y}_1,\ldots,\mathbf
{y}_L)$ under the model corresponding to a particular value $\vartheta$ is given by
%
%
\begin{equation}
\label{eq-GL}
\mathscr{L}_{\overline{\mathbf{y}}^L}(\vartheta) = (2\uppi)^{-L d/2}
\Biggl(\prod
_{n=1}^L{\det V_{\vartheta,n}} \Biggr)^{-1/2}\exp\Biggl\{-\frac{1}{2}\sum
_{n=1}^L{\be_{\vartheta,n}^TV_{\vartheta,n}^{-1}\be_{\vartheta,n}}
\Biggr\},
\end{equation}
where $\be_{\vartheta,n}$ is the residual of the minimum mean-squared
error linear predictor of $\mathbf{y}_n$ given the preceding
observations, and
$V_{\vartheta,n}$ is the corresponding covariance matrix. From a
practical perspective, it is important to note that all quantities
necessary to evaluate the Gaussian likelihood $\eqref{eq-GL}$ can be
conveniently computed by using the Kalman recursions (\cite{brockwell1991tst}, Section 12.2) and the state space representation given in
Lemma \ref{sampledSSM}. In case the observations $\overline{\mathbf
{y}}^L$ are
(part of) a realization of the sampled MCARMA process $\Y
^{(h)}_{\vartheta_0}$ corresponding to the parameter value $\vartheta
_0$, the prediction error sequence $ (\be_{\vartheta_0,n}
)_{n\geq1}$ is -- up to an additive, exponentially decaying term which
comes from the initialization of the Kalman filter -- (part of) a
realization of the innovations sequence $\bepsh$ of $\Y
^{(h)}_{\vartheta
_0}$. In order to be able to analyze the
asymptotic behavior of the natural GML estimator
\[
\hat\vartheta^L = \operatorname{argmax}\limits_{\vartheta\in\Theta
}\mathscr
{L}_{\overline{\mathbf{y}}^L}(\vartheta)
\]
in the limit as $L\to\infty$, it is necessary to have a central limit
theorem for sums of the~form
%
%
\begin{equation}
\label{eq-CLT}
\frac{1}{\sqrt{L}}\sum_{n=1}^L{\frac{\partial}{\partial\vartheta
} [\log\det V_{\vartheta,n}+\be_{\vartheta,n}^TV_{\vartheta
,n}^{-1}\be_{\vartheta,n} ] \bigg|_{\vartheta=\vartheta_0}}.
\end{equation}
Existing results in the literature \cite
{bradley2007introduction,herrndorf1984} ensure that various notions of
weak dependence, and, in particular, strong mixing, are sufficient for
a central limit theorem for the expression \eqref{eq-CLT} to hold.
Theorem \ref{innoMix} below is thus the necessary starting point for
the development of an estimation theory for multivariate CARMA
processes which involves some additional issues like identifiability of
parametrizations and is thus beyond the scope of this paper.

Before presenting the sufficient condition for the innovations $\bepsh$
to be completely regular, we first observe that the eigenvalues
$\lambda
_1,\ldots,\lambda_N$ of $A$ are the roots of the characteristic
polynomial $z\mapsto\det(z\I_N-A)$, which, by the fundamental theorem
of algebra, implies that they are either real\vadjust{\goodbreak} or occur in complex
conjugate pairs. We can therefore assume that they are ordered in such
a way that for some $r\in\{0,\ldots,N\}$,
\[
\lambda_\nu\in\R,\ 1\leq\nu\leq r, \qquad\lambda_\nu
=\overline
{\lambda_{\nu+1}}\in\C\backslash\R,\qquad\nu=r+1,r+3,\ldots,N-1.
\]
By Lebesgue's decomposition theorem \cite{klenke2008probability}, Theorem 7.33, every measure $\mu$ on $\R^d$ can be
uniquely decomposed\vadjust{\goodbreak} as $\mu=\mu_{\text{c}}+\mu_{\text{s}}$, where
$\mu
_{\text{c}}$ and $\mu_{\text{s}}$ are absolutely continuous and
singular, respectively, with respect to the $d$-dimensional Lebesgue
measure. If $\mu_{\text{c}}$ is not the zero measure, then we say that
$\mu$ has a non-trivial absolutely continuous component.

\begin{theorem}
\label{innoMix}
Assume that $\mathbf{Y}$ is the output process of the continuous-time state
space model $(A,B,C,\Lb)$ satisfying Assumptions  \textup{\ref{assumlevy2mom}}, \textup{\ref{assumEigen1}} and \textup{\ref{assumEigen2}}. Denote by $\bepsh$ the innovations
of the sampled process $\Y^{(h)}$ and further assume that the law of
the $\R^{mN}$-valued random variable
%
%
\begin{equation}
\label{defM}
\vecM^{(h)}= \left[
\matrix{{\M^{(h)}_1}^{\mathrm{T}} & \cdots& {\M^{(h)}_r}^{\mathrm{T}} &
{\underline\M^{(h)}_{r+1}}^{\mathrm{T}} & {\underline\M^{(h)}_{r+3}}^{\mathrm{T}} &
\cdots&
{\underline\M^{(h)}_{N-1}}^{\mathrm{T}}
}
\right]^{\mathrm{T}},
\end{equation}
where
%
%
\begin{equation}
\hspace*{-15pt}\underline\M^{(h)}_\nu= \left[
\matrix{{\re\M^{(h)}_\nu}^{\mathrm{T}} & {\imag\M^{(h)}_\nu}^{\mathrm{T}}
}
\right]^{\mathrm{T}},\quad\M^{(h)}_\nu=\int_0^h{\mathrm{e}^{(h-u)\lambda
_\nu}\,\mathrm{d}\Lb
(u)},\quad\nu=1,\ldots,N,
\end{equation}
has a non-trivial absolutely continuous component with respect to the
$mN$-dimensional Lebesgue measure. Then, $\bepsh$ is exponentially
completely regular.
\end{theorem}

The assumption on the distribution of $\vecM^{(h)}$ made in Theorem \ref
{innoMix} is not very restrictive. Its verification is based on the
following lemma, which allows us to derive sufficient conditions in
terms of the L\'{e}vy process $\Lb$ which show that it is indeed
satisfied in most practical situations.
\begin{lemma}
\label{MOU}
There exist matrices $G\in M_{mN}(\R)$ and $H\in M_{mN,m}(\R)$ such
that $\vecM^{(h)}=\vecM(h)$, where $(\vecM(t))_{t\geq0}$ is the unique
solution to the stochastic differential equation
%
%
\begin{equation}
\label{MOU-eq}
\mathrm{d}\vecM(t)=G\vecM(t)\,\mathrm{d}t+H\,\mathrm{d}\Lb
(t),\qquad\vecM
(0)=\bzero_{mN}.
\end{equation}
Moreover, $\rank H = m$ and the $mN\times mN$ matrix $ \left[
H \enskip GH \enskip\cdots\enskip G^{N-1}H
\right]$ is non-singular.
\end{lemma}

The last part of the statement is referred to as \textit{controllability}
of the pair $(G,H)$ and is essential in the proofs of the following
explicit sufficient conditions for Theorem \ref{innoMix} to hold.
\begin{proposition}
\label{GaussianCond}
Assume that the L\'{e}vy process $\Lb$ has a non-singular Gaussian
covariance matrix $\Sigma^{\mathcal{G}}$. Theorem \ref{innoMix} then holds.
\end{proposition}

\begin{pf}
By \cite{sato2006additive}, Corollary 2.19, the law of $\vecM^{(h)}$ is
infinitely divisible with Gaussian covariance matrix given by $\int
_0^h{\mathrm{e}^{Gu}H\Sigma^{\mathcal{G}}H^{\mathrm{T}}\mathrm{e}^{G^\mathrm{T}u}\,
\mathrm{d}u}$. By the
controllability of $(G,H)$ and~\cite{bernstein2005matrix}, Lemma~12.6.2,
this matrix is non-singular and \cite{sato1999lpa}, Exercise 29.14 completes the proof.
\end{pf}

A simple L\'{e}vy process of practical importance which does not have a
non-singular Gaussian covariance matrix is the \textit{compound Poisson
Process}, which is defined by $\Lb(t)=\sum_{n=1}^{N(t)}{\J_n}$, where
$(N(t))_{t\in\R^+}$ is a Poisson process and $(\J_n)_{n\in\Z}$ is an
i.i.d. sequence independent of $(N(t))_{t\in\R^+}$; the law of $\J_n$
is called the \textit{jump size distribution}. The proof of \cite{priola2009densities}, Theorem 1.1, in conjunction with Lemma \ref
{MOU}, implies the following result.\vadjust{\goodbreak}
\begin{proposition}
Assume that $\Lb$ is a compound Poisson process with absolutely
continuous jump size distribution. Theorem \ref{innoMix} then holds.
\end{proposition}

Under a similar smoothness assumption, the conclusion of Theorem \ref
{innoMix} also holds in the case of infinite activity L\'{e}vy
processes. The statement follows from applying \cite{priola2009densities}, Theorem 1.1, to equation \eqref{MOU-eq}.
\begin{proposition}
\label{DensityCond}
Assume that the L\'{e}vy measure $\nuL$ of $\Lb$ satisfies $\nuL(\R
^m)=\infty$ and that there exists a positive constant $\rho$ such that
$\nuL$ restricted to the ball $\{\bx\in\R^m\dvtx\|\bx\|\leq\rho
\}$ has a
density with respect to the $m$-dimensional Lebesgue measure. Theorem
\ref{innoMix} then holds.\looseness=-1
\end{proposition}

While the preceding three propositions already cover a wide range of L\'
{e}vy processes encountered in practice, there are some relevant cases
which are not yet taken care of, in particular, the construction of the
L\'{e}vy process as a vector of independent univariate L\'{e}vy
processes (Corollary~\ref{indL} below).
To also cover this and related choices, we employ the polar
decomposition for L\'{e}vy measures \cite{barndorff2006some}, Lemma 2.1. By this result, for every L\'{e}vy measure
$\nuL$, there exists a probability measure $\alpha$ on the
$(m-1)$-sphere $S^{m-1}\coloneqq\{\bx\in\R^m\dvtx\|\bx\|=1\}$ and
a family
$\{\nu_{\bxi}\dvtx{\bxi}\in S^{m-1}\}$ of measures on $\R^+$ such
that for
each Borel set $B\in\mathcal{B}(\R^+)$, the function ${\bxi}\mapsto
\nu
_{\bxi}(B)$ is measurable and
%
%
\begin{equation}
\label{polardecomp}
\nuL(B)=\int_{S^{m-1}}{\int_0^\infty{I_B(\lambda{\bxi})\nu_{\bxi
}(\mathrm{d}\lambda)}\alpha(\mathrm{d}{\bxi})},\qquad B\in\mathcal
{B}(\R
^m\backslash\{\bzero_m\}).
\end{equation}
A hyperplane in a finite-dimensional vector space is a linear subspace
of codimension~one.

\begin{proposition}
If the L\'{e}vy measure $\nuL$ has a polar decomposition $(\alpha,\nu
_{\bxi}\dvtx\bxi\in S^{m-1})$ such that for any hyperplane $\mathcal
{H}\subset\R^m$, it holds that $\int_{S^{m-1}}{I_{\R^m\backslash
\mathcal
{H}}({\bxi})\int_0^\infty{\nu_{\bxi}(\mathrm{d}\lambda)}\alpha
(\mathrm
{d}{\bxi})}=\infty$, then Theorem~\ref{innoMix} holds.
\end{proposition}

\begin{pf}
The proof rests on the main theorem of \cite{simon2010absolute}. We
denote by $\operatorname{im} H$ the image of the linear operator
associated with the matrix $H$. Since $\rank H=m$ and the pair $(G,H)$
is controllable, we only have to show that $\nuL(\{\bx\in\R^m\dvtx
H\bx\in
\operatorname{im}{H}\backslash\mathcal{H}\})=\infty$ for all
hyperplanes $\mathcal{H}\subset\operatorname{im}{H}$, and since $\R
^m\cong\operatorname{im}{H}$, the last condition is equivalent to
$\nuL
(\R^m\backslash\mathcal{H})=\infty$ for all hyperplanes $\mathcal
{H}\subset\R^m$. Using equation \eqref{polardecomp} and the fact that for every
${\bxi}\in S^{m-1}$ and every $\lambda\in\R^+$, the vector $\lambda
{\bxi}$ is in $\mathcal{H}$ if and only if the vector ${\bxi}$ is, this
is seen to be equivalent to the assumption of the proposition.
\end{pf}

\begin{corollary}
If the L\'{e}vy measure $\nuL$ has a polar decomposition $(\alpha,\nu
_{\bxi}\dvtx\bxi\in S^{m-1})$ such that $\alpha(S^{m-1}\backslash
\mathcal
{H})$ is positive for all hyperplanes $\mathcal{H}\in\R^m$ and $\nu
_{\bxi}(\R^+)=\infty$ for $\alpha$-almost every $\bxi$, then Theorem
\ref{innoMix} holds.
\end{corollary}

\begin{corollary}
If the L\'{e}vy measure $\nuL$ has a polar decomposition $(\alpha,\nu
_{\bxi}\dvtx\bxi\in S^{m-1})$ such that for some linearly independent
vectors ${\bxi}_1,\ldots,{\bxi}_m\in S^{m-1}$, it holds that $\alpha
({\bxi}_k)>0$ and $\nu_{{\bxi}_k}(\R^+)=\infty$ for $k=1,\ldots,m$,
then Theorem \ref{innoMix} holds.\vadjust{\goodbreak}
\end{corollary}

\begin{corollary}
\label{indL}
Assume that $l\geq m$ is an integer and that the matrix $R\in
M_{m,l}(\R
)$ has full rank~$m$. If $\Lb=R \left(
L_1 \enskip\cdots\enskip L_l
\right)^{\mathrm{T}}$, where $L_k$, $k=1,\ldots, l$, are independent univariate
L\'
{e}vy processes with L\'{e}vy measures $\nu^L_k$ satisfying $\nu
^L_k(\R
)=\infty$, then Theorem \ref{innoMix} holds.\vspace*{-3pt}
\end{corollary}

%
\section{Proofs}\label{sec5}\vspace*{-3pt}
\subsection{\texorpdfstring{Proofs for Section \protect\ref{section-MCARMA}}{Proofs for Section
3}}\vspace*{-3pt}
\label{section-proofs}

\begin{pf*}{Proof of Theorem \ref{equivalenceMCARMASSM}}
The first step of the proof is to show that any pair $(P,Q)$ of the
form \eqref{ARMApoly} is a left matrix fraction description of
$\mathcal
{C}(z\I_{pd}-\mathcal{A})^{-1}\mathcal{B}$, provided $\mathcal{A}$,
$\mathcal{B}$ and~$\mathcal{C}$ are defined as in equations \eqref
{MCARMAcoeffABC}. We first show the relation
%
%
\begin{equation}
\label{MFDeq}
(z \I_{pd}-\mathcal{A})^{-1}\mathcal{B} = \left[
\matrix{w_1(z)^{\mathrm{T}} & \cdots& w_p^{\mathrm{T}}(z)
}
\right]^{\mathrm{T}},
\end{equation}
where $w_j(z)\in M_{d,m}(\R\{z\})$, $j=1,\ldots, p$, are defined by
the equations
%
%
\begin{subequation}
\label{MFDdefw}
\begin{equation}
w_j(z)=\frac{1}{z}\bigl(w_{j+1}(z)+\beta_j\bigr),\qquad j=1,\ldots,p-1,
\end{equation}
and\vspace*{-1pt}
\begin{equation}
w_p(z)=\frac{1}{z} \Biggl(-\sum_{k=0}^{p-1}{A_{p-k}w_{k+1}(z)}+\beta
_p \Biggr).
\end{equation}
\end{subequation}
Since it has been shown in \cite{marquardt2007multivariate}, Theorem 3.12, that $w_1(z)=P(z)^{-1}Q(z)$ this will
prove the assertion. Equation \eqref{MFDeq} is clearly equivalent to
$\mathcal{B}=(z\I_{pd}-\mathcal{A}) \left[
w_1(z)^{\mathrm{T}} \enskip\cdots\enskip w_p^{\mathrm{T}}(z)
\right]^{\mathrm{T}}$, which explicitly reads
\begin{eqnarray*}
\beta_j &=& z w_j(z)-w_{j+1}(z),\qquad j=1,\ldots, p-1,\\
\beta_p &= &z w_p(z)+A_p w_1(z)+\cdots+ A_1 w_p(z)
\end{eqnarray*}
and is thus equivalent to equations \eqref{MFDdefw}.

For the second step consider a given state space model $(A,B,C,\Lb)$.
Using the spectral representation \cite{lax2002functional}, Theorem 17.5,
%
%
\begin{equation}
\label{specrep}
\mathrm{e}^{At}=\frac{1}{2\uppi \mathrm{i}}\int_\Gamma{\mathrm{e}^{zt}(z\I
_N-A)^{-1}\,\mathrm
{d}z},\qquad t\in\R,
\end{equation}
where $\Gamma$ is some closed contour in $\C$ winding around each
eigenvalue of $A$ exactly once, it follows that\vspace*{-1pt}
\begin{eqnarray*}
\Y(t)&=&\int_{-\infty}^t{C\mathrm{e}^{A(t-u)}B\,\mathrm{d}\Lb
(u)}=\frac{1}{2\uppi\mathrm{i}}\int_{-\infty}^t{\int_\Gamma{\mathrm{e}^{z(t-u)}C(z\I
_N-A)^{-1}B\,\mathrm
{d}z}\,\mathrm{d}\Lb(u)}\\
&=&\frac{1}{2\uppi\mathrm{i}}\int_{-\infty}^t{\int_\Gamma
{\mathrm{e}^{z(t-u)}P(z)^{-1}Q(z)\,\mathrm{d}z}\,\mathrm{d}\Lb
(u)}\\
&=&\frac{1}{2\uppi\mathrm{i}}\int_{-\infty}^t{\int_\Gamma{\mathrm{e}^{z(t-u)}\mathcal{C}(z\I
_{pd}-\mathcal
{A})^{-1}\mathcal{B}\,\mathrm{d}z}\,\mathrm{d}\Lb(u)}\\
&=&\int_{-\infty}^t{\mathcal{C}\mathrm{e}^{\mathcal
{A}(t-u)}\mathcal{B}\,\mathrm
{d}\Lb(u)},\
\end{eqnarray*}
where $\mathcal{A}$, $\mathcal{B}$ and $\mathcal{C}$ are defined in
terms of $(P,Q)$ by equations \eqref{MCARMAcoeffABC}. Thus $\Y$ is an MCARMA
process with autoregressive polynomial $P$ and moving average
polynomial~$Q$.\vspace*{3pt}
\end{pf*}

%
\subsection{\texorpdfstring{Proofs for Section \protect\ref{section-results}}{Proofs for Section 4}}
\vspace*{3pt}

In this section we present the proofs of our main results, Theorem \ref
{weakARMAthm}, Theorem \ref{innoMix} and Lemma~\ref{MOU}, as well as
several auxiliary results. The first is a generalization of \cite{brockwell2010estimation}, Proposition 2, expressing MCARMA processes as
a sum of multivariate Ornstein--Uhlenbeck processes.
\begin{proposition}
\label{decompositionCT}
Let $\Y$ be the the output process of the state space system \eqref
{ssm} and assume that Assumption~\textup{\ref{assumEigen2}} holds. Then, there exist
vectors $\bs_1,\ldots,\bs_N\in\C^m\backslash\{\bzero_m\}$ and
$\bb
_1,\ldots,\bb_N\in\C^d\backslash\{\bzero_d\}$ such that $\Y$ can be
decomposed into a sum of dependent, complex-valued Ornstein--Uhlenbeck
processes as $ \Y(t)=\sum_{\nu=1}^N{\Y_\nu(t)}$, where
%
%
\begin{equation}
\label{decomprec}
\Y_\nu(t)=\mathrm{e}^{\lambda_\nu(t-s)}\Y_\nu(s)+\bb_\nu\int
_{s}^t{ \mathrm{e}^{\lambda
_\nu(t-u)}\,\mathrm{d}\langle\bs_\nu,\Lb(u)\rangle},\qquad s,t\in
\R,\  s<t.
\end{equation}
\end{proposition}

\begin{pf}
We first choose a left matrix fraction description $(P,Q)$ of the
transfer function $z\mapsto C(z\I_N-A)^{-1}B$ such that $z\mapsto\det
P(z)$ and $z\mapsto\det Q(z)$ have no common zeros and $z\mapsto\det
P(z)$ has no multiple zeros. This is always possible, by Assumption
\ref
{assumEigen2}. Inserting the spectral representation \eqref{specrep} of
$\mathrm{e}^{At}$ into the kernel $g(t)$ (equation \eqref{ssmMA}), we get
$g(t)=\frac{1}{2\uppi\mathrm{i}}\int_\Gamma{\mathrm{e}^{zt}C(z\I_N-A)^{-1}B\,
\mathrm
{d}z}I_{[0,\infty)}(t)$ and, by construction, the integrand equals
$\mathrm{e}^{zt}P(z)^{-1}Q(z)I_{[0,\infty)}(t)$. After writing
$P(z)^{-1}=\frac
{1}{\det P(z)}\adj P(z)$, where $\adj$ denotes the adjugate of a
matrix, an elementwise application of the residue theorem from complex
analysis (\cite{dieudonne2006foundations},\linebreak[2]  Theorem 9.16.1) shows that
\[
g(t)=\sum_{\nu=1}^N{\mathrm{e}^{\lambda_\nu t}\frac{1}{(\det
P)'(\lambda_\nu
)}\adj P(\lambda_\nu)Q(\lambda_\nu)}I_{[0,\infty)}(t),
\]
where $(\det P)'(\lambda_\nu)\coloneqq\frac{\mathrm{d}}{\mathrm{d}z}\det
P(z)|_{z=\lambda
_\nu}$ is non-zero because $z\mapsto\det P(z)$ has only simple zeros.
The same fact, in conjunction with the Smith decomposition of $P$
(\cite{bernstein2005matrix}, Theorem 4.7.5), also implies that $\rank
P(\lambda_\nu)=d-1$ and thus $\rank\adj P(\lambda_\nu)=1$ (\cite{bernstein2005matrix}, Fact 2.14.7(ii)). Since $\det P$ and
$\det Q$
have no common zeros, $[(\det P)'(\lambda_\nu)]^{-1}\adj P(\lambda
_\nu
)Q(\lambda_\nu)$ also has rank one and can thus be written as $\bb
_\nu
\bs_\nu^{\mathrm{T}}$ for some non-zero $\bs_\nu\in\C^m$ and $\bb_\nu\in
\C^d$
(\cite{halmos1974finite}, Section 51, Theorem~1).
\noqed\mbox{}\qed
\end{pf}

\begin{lemma}
\label{sampledSSM}
Assume that $\Y$ is the output process of the state space model \eqref
{ssm}. The sampled process $\Y^{(h)}$ then has the state space representation
%
%
\begin{equation}
\label{sampledSSM-eq}
\X_n = \mathrm{e}^{Ah}\X_{n-1}+\NN_n,\qquad\NN_n= \int
_{(n-1)h}^{nh}{\mathrm{e}^{A(nh-u)}B\,\mathrm{d}\Lb_u},\qquad\Y
^{(h)}_n = C\X^{(h)}_n.
\end{equation}
The sequence $(\NN_n)_{n\in\Z}$ is i.i.d. with mean zero and
covariance matrix
%
%
\begin{equation}
\label{discnoisecov}
\cancel\Sigma=\E\NN_n\NN_n^{\mathrm{T}}=\int_0^h{\mathrm{e}^{Au}B\Sigma
^{\Lb}
B^{\mathrm{T}}\mathrm{e}^{A^\mathrm{T}u}\,\mathrm{d}u}.
\end{equation}
\end{lemma}

\begin{pf}
Equations \eqref{sampledSSM-eq} follow from setting $t=nh$, $s=(n-1)h$
in equation \eqref{markovstate}. It is an immediate consequence of the
L\'{e}vy process $\Lb$ having independent, homogeneous increments that
the sequence $(\mathbf{N}_n)_{n\in\Z}$ is i.i.d. and that its covariance
matrix $\cancel\Sigma$ is given by equation \eqref{discnoisecov}.
\end{pf}

From this, we can now proceed to prove the weak vector ARMA
representation of the process~$\Y^{(h)}$.
\begin{pf*}{Proof of Theorem \ref{weakARMAthm}}
It follows from setting $t=nh$, $s=(n-1)h$ in equation \eqref{decomprec} that
$\Y
_n^{(h)}$ can be decomposed as $\Y_n^{(h)}=\sum_{\nu=1}^N{\Y_{\nu
,n}^{(h)}}$, where $\Y_\nu^{(h)}$, satisfying
\[
\label{discrec}
\Y_{\nu,n}^{(h)} = \mathrm{e}^{\lambda_\nu h}\Y_{\nu
,n-1}^{(h)}+\ZZ_{\nu
,n}^{(h)},\qquad\ZZ_{\nu,n}^{(h)}=\bb_\nu\int
_{(n-1)h}^{nh}{\mathrm{e}^{\lambda
_\nu(nh-u)}\,\mathrm{d}\langle\bs_\nu,\Lb(u)\rangle},
\]
are the sampled versions of the component $\operatorname{MCAR}(1)$
processes from Proposition \ref{decompositionCT}. Analogously to \cite{brockwell2009existence}, Lemma 2.1, we can show by induction that for
each $k\in\N_0$ and all complex $d\times d$ matrices $c_1,\ldots,c_k$,
it holds that
%
%
\begin{eqnarray}
\label{inductionhyp}
\Y_{\nu,n}^{(h)}&=&\sum_{r=1}^k{c_r\Y_{\nu,n-r}^{(h)}}+ \Biggl[\mathrm
{e}^{\lambda
_\nu hk}-\sum_{r=1}^k{c_r\mathrm{e}^{\lambda_\nu h(k-r)}} \Biggr]\Y_{\nu
,n-k}^{(h)}\nonumber
\\[-8pt]
\\[-8pt]
&&{}+\sum_{r=0}^{k-1}\Biggl [\mathrm{e}^{\lambda_\nu hr}-\sum
_{j=1}^r{c_j\mathrm{e}^{\lambda_\nu h(r-j)}} \Biggr]\ZZ_{\nu,n-r}^{(h)}.
\nonumber
\end{eqnarray}
If we then use the fact that $\mathrm{e}^{-h\lambda_\nu}$ is a root
of $z\mapsto
\varphi(z)$, which means that $\mathrm{e}^{Nh\lambda_\nu}-\varphi
_1\mathrm{e}^{(N-1)h\lambda_nu}-\cdots-\varphi_N = 0$, and set
$k=N$, $c_r=\I
_d\varphi_r$, then equation \eqref{inductionhyp} becomes
\[
\varphi(B)\Y_{\nu,n}^{(h)}=\sum_{r=0}^{N-1} \Biggl [\mathrm
{e}^{rh\lambda_\nu}-\sum
_{j=1}^r{\varphi_j\mathrm{e}^{\lambda_\nu h(r-j)}} \Biggr]\ZZ_{\nu,n-r}^{(h)}.
\]
Summing over $\nu$ and rearranging shows that this can be written as
%
%
\begin{equation}
\label{preMA}
\varphi(B)\Y_n^{(h)}=\sum_{\nu=1}^N\V_{\nu,n-\nu+1}^{(h)},
\end{equation}
where the i.i.d. sequences $ (\V_{\nu,n}^{(h)} )_{n\in\mathbb
{Z}}$, $\nu\in\{1,\ldots,N\}$, are defined by
%
%
\begin{equation}
\label{DefV}
\V_{\nu,n}^{(h)}=\int_{(n-1)h}^{nh}\sum_{\mu=1}^N{\bb_\mu
\Biggl[\mathrm{e}^{\lambda_\mu h(\nu-1)}-\sum_{\kappa=1}^{\nu
-1}{\varphi_\kappa
\mathrm{e}^{\lambda_\mu h(\nu-\kappa-1)}} \Biggr]\mathrm{e}^{\lambda
_\mu(nh-u)}}\,\mathrm
{d}\langle\bs_\mu,\Lb(u)\rangle.
\end{equation}
By a straightforward generalization of \cite{brockwell1991tst}, Proposition 3.2.1, there exists a monic Schur-stable polynomial
$\Theta(z)=\mathrm{I}_d+\Theta_1z+\cdots+\Theta_{N-1}z^{N-1}$ and a
white noise sequence $\tilde{\beps}$ such that the $(N-1)$-dependent
sequence $\varphi(B)\Y^{(h)}$ has the moving average representation
$\varphi(B)\Y_n^{(h)}=\Theta(B)\tilde{\beps}_n$. Since both
$\varphi$
and $\Theta$ are monic, and $\varphi$ is Schur stable (by Assumption
\ref{assumEigen1}), $\tilde{\beps}$ is the innovation process of
$\Y
^{(h)}$ and so it follows that $\tilde{\beps}=\bepsh$ because the
innovations of a stochastic process are uniquely determined.
\end{pf*}

As a corollary, we obtain that the innovations sequence $\bepsh$ itself
satisfies a set of strong VARMA equations, the attribute \textit{strong}
referring to the fact that the noise sequence is i.i.d., not merely
white noise.
\begin{corollary}
\label{InnoARMA-thm}
Assume that $\Y$ is the output process of the state space system
\textup{\eqref{ssm}} satisfying Assumptions \textup{\ref{assumlevy2mom}}, \textup{\ref{assumEigen1}} and \textup{\ref{assumEigen2}}. Further assume that $\bepsh$ is the innovations sequence
of the sampled process $\Y^{(h)}$. There then exists a monic,
Schur-stable polynomial $\Theta\in M_d(\C[z])$ of degree at most $N-1$,
a polynomial $\theta\in M_{d,dN}(\R[z])$ of degree $N-1$ and a
$C^{dN}$-valued i.i.d. sequence $\W^{(h)}=(\W_n^{(h)})_{n\in\Z}$,
such that
%
%
\begin{equation}
\label{innoARMA}
\Theta(B)\bepsh_n=\theta(B)\W_n^{(h)},\qquad n\in\Z.
\end{equation}
\end{corollary}

\begin{pf}
Combining equations \eqref{weakARMA} and   \eqref{preMA} gives
%
%
\begin{eqnarray}
\label{preinnoARMA}
&&\bepsh_n+\Theta^{(h)}_1\beps_{n-1}+\cdots+\Theta^{(h)}_{N-1}\beps
_{n-N+1}\nonumber
\\[-8pt]
\\[-8pt] && \quad = \V_{1,n}^{(h)}+\V_{2,n-1}^{(h)}+\cdots+\V
_{N,n-N+1}^{(h)},\qquad n\in\Z,
\nonumber
\end{eqnarray}
and with the definitions
%
%
\begin{subequation}
\begin{eqnarray}
\label{DefW}\W_n^{(h)}&=&\left[
\matrix{{\V_{1,n}^{(h)}}^{\mathrm{T}} & \cdots& {\V_{N,n}^{(h)}}^{\mathrm{T}}
}
\right]^{\mathrm{T}}\in\C^{dN},\qquad n\in\Z,\\
\label{Deftheta}\theta(z)&=&\sum_{j=1}^N{\theta_jz^{j-1}},\nonumber
\\[-8pt]
\\[-8pt]
\theta_\nu
&=&[\vphantom{\I_d} \underbrace{
\matrix{0_d & \cdots& 0_d
}
}_{\nu-1\mathrm{\ times}}
\;\;\; \I_d\;\;\;
\underbrace{
\matrix{ 0_d & \cdots& 0_d
}
}_{N-\nu\mathrm{\ times}} ]\in M_{d,dN}(\R),\qquad\nu=1,\ldots,N,
\nonumber
\end{eqnarray}
\end{subequation}
equation \eqref{preinnoARMA} becomes $\Theta(B)\bepsh_n=\theta(B)\W
_n^{(h)}$, showing that $\bepsh$ is indeed a vector ARMA process.
\end{pf}

This corollary is the central step in establishing complete regularity
of the innovations process~$\bepsh$.\vadjust{\goodbreak}
\begin{pf*}{Proof of Theorem \ref{innoMix}}
We define the $\R^{mN}$-valued random variables
\[
\vecM^{(h)}_n= \left[
\matrix{{\M^{(h)}_{n,1}}^{\mathrm{T}} & \cdots& {\M^{(h)}_{n,r}}^{\mathrm{T}}
& {\underline\M^{(h)}_{n,r+1}}^{\mathrm{T}} & {\underline\M^{(h)}_{n,r+3}}^{\mathrm{T}} &
\cdots& {\underline\M^{(h)}_{n,N-1}}^{\mathrm{T}}
}
\right]^{\mathrm{T}},\qquad n\in\Z,
\]
where
\begin{eqnarray*}
\underline\M^{(h)}_{n,\nu} &=& \left[
\matrix{{\re\M^{(h)}_{n,\nu}}^{\mathrm{T}} & {\imag\M^{(h)}_{n,\nu}}^{\mathrm{T}}
}
\right]^{\mathrm{T}},\\
\M^{(h)}_{n,\nu}&=&\int_{(n-1)h}^{nh}{\mathrm
{e}^{\lambda_\nu
(nh-u)}\,\mathrm{d}\Lb(u)},\qquad\nu=1,\ldots,N,\  n\in\Z.
\end{eqnarray*}
Clearly, the sequence $(\vecM^{(h)}_n)_{n\in\Z}$ is i.i.d. and
$\vecM
^{(h)}$ is equal to $\vecM^{(h)}_1$. We now argue that the vector $\W
^{(h)}_n$, as defined in equation \eqref{DefW}, is equal to a linear
transformation of $\vecM^{(h)}_n$. By equation~\eqref{DefV}, $\W
^{(h)}_n=[\Gamma^{\mathrm{T}}\otimes\I_d] \left[
(\bb_1\bs_1^{\mathrm{T}} \M^{(h)}_{n,1})^{\mathrm{T}} \enskip\cdots\enskip(\bb_N\bs
_N^{\mathrm{T}} \M^{(h)}_{n,N})^{\mathrm{T}}
\right]^{\mathrm{T}}$, where $\Gamma=(\gamma_{\mu,\nu})\in M_N(\C)$ is given by
$\gamma_{\mu,\nu}=\mathrm{e}^{\lambda_\mu h(\nu-1)}+\sum_{\kappa
=1}^{\nu
-1}{\varphi_\kappa\mathrm{e}^{\lambda_\mu h(\nu-\kappa-1)}}$.
With the notation
\[
B= \pmatrix{\displaystyle
\bb_1 & \bzero_d & \ldots& \bzero_d\cr\displaystyle
\bzero_d & \bb_2 & \ddots& \vdots\cr\displaystyle
\vdots&\ddots&\ddots& \bzero_d \cr\displaystyle
\bzero_d & \ldots& \bzero_d & \bb_N
}
\in M_{dN,N}(\C),\qquad S=\pmatrix{\displaystyle
\bs_1^{\mathrm{T}} & \bzero_d^{\mathrm{T}} & \ldots& \bzero_d^{\mathrm{T}}\cr\displaystyle
\bzero_d^{\mathrm{T}} & \bs_2^{\mathrm{T}} & \ddots& \vdots\cr\displaystyle
\vdots&\ddots&\ddots& \bzero_d^{\mathrm{T}} \cr\displaystyle
\bzero_d^{\mathrm{T}} & \ldots& \bzero_d^{\mathrm{T}} & \bs_N^{\mathrm{T}}
}
\in M_{N,mN}(\C),
\]
we get $ \left[
(\bb_1\bs_1^{\mathrm{T}} \M^{(h)}_{n,1} )^{\mathrm{T}} \enskip\cdots\enskip
(\bb_N\bs_N^{\mathrm{T}} \M^{(h)}_{n,N} )^{\mathrm{T}}
\right]^{\mathrm{T}}=BS \left[
{\M^{(h)}_{n,1}}^{\mathrm{T}} \enskip\cdots\enskip{\M^{(h)}_{n,N}}^{\mathrm{\!\!T}}
\right]^{\mathrm{T}}$.
We recall that for
$\nu = r + 1, r + 3,\ldots,N - 1$, the eigenvalues of $A$ satisfy $\lambda_\nu = \overline{\lambda_{\nu+1}} \in
\mathbb{C}\backslash \mathbb{R}$, which implies
that
\[
\mathbf{M}^{(h)}_{n,\nu}=\operatorname{Re}\mathbf{M}^{(h)}_{n,\nu}+\mathrm{i}\operatorname{Im}\mathbf{M}^{(h)}_{n,\nu}
 \quad \mbox{and} \quad \mathbf{M}^{(h)}_{n,\nu+1}=\overline{\mathbf{M}^{(h)}_{n,\nu}} =\operatorname{Re}\mathbf{M}^{(h)}_{n,\nu}
 -\mathrm{i}\operatorname{Im}\mathbf{M}^{(h)}_{n,\nu}.
\]
Consequently, we obtain that
$ \left[
{\M^{(h)}_{n,1}}^{\mathrm{T}} \enskip\cdots\enskip{\M^{(h)}_{n,N}}^{\!\!\mathrm{T}}
\right]^{\mathrm{T}}=[K\otimes\I_m]\vecM^{(h)}_n$, where
\[
K= \pmatrix{\displaystyle
\I_r & & & \cr\displaystyle
&J & & \cr\displaystyle
& &\ddots& \cr\displaystyle
& & &J
}
\in M_N(\C), \qquad J= \pmatrix{\displaystyle1&\mathrm{i}\cr\displaystyle1&-\mathrm{i}
}
,
\]
so that, in total, $\W^{(h)}_n= F\vecM^{(h)}_n$ with $ F=[\Gamma
^{\mathrm{T}}\otimes\I_d]BS[K\otimes\I_m]\in M_{dN,mN}(\C)$. It follows that the
VARMA equation \eqref{innoARMA} for $\bepsh$ becomes $\Theta
(B)\bepsh_n
= \tilde\theta(B)\vecM^{(h)}_n$, where $\tilde\theta(z)= \theta(z)F$.
By the invertibility of $\Theta$, the transfer function $k\dvtx
z\mapsto
\Theta(z)^{-1}\tilde\theta(z)$ is analytic in a disc containing the
unit disc and permits a power series expansion $k(z)=\sum_{j=0}^\infty
{\Psi_j z^j}$. We next argue that the impulse responses $\Psi_j$ are
necessarily \textit{real} $d\times mN$ matrices. Since both $\bepsh
_n$ and
$\vecM^{(h)}_n$ are real-valued, it follows from taking the imaginary
part of the equation $\bepsh_n=k(B)\vecM^{(h)}_n$ that $\bzero
_d=\sum
_{j=0}^\infty{\imag\Psi_j\vecM^{(h)}_{n-j}}$. Consequently, $0=\Cov
(\bzero_d)=\sum_{j=0}^\infty{\imag\Psi_j\Cov(\vecM
^{(h)}_{n-j})\imag\Psi
_j^{\mathrm{T}}}$ and since each term in the sum is a positive semidefinite
matrix, it follows that
$\imag\Psi_j\Cov(\vecM^{(h)}_{n-j})\imag\Psi_j^{\mathrm{T}}=0$ for every
$j$. The
existence of an absolutely continuous component of the law of $\vecM
^{(h)}_{n-j}$ with respect to the $mN$-dimensional Lebesgue measure
implies that $\Cov(\vecM^{(h)}_{n-j})$ is non-singular and it thus
follows that $\imag\Psi_j=0$ for every $j$. Hence, $k(z)\in
M_{d,mN}(\R
)$ for all real $z$, and consequently $k\in M_{d,mN}(\R\{z\})$. \cite{hannan1987stl}, Theorem 1.2.1(iii), then implies that there
exists a~stable $(\vecM^{(h)}_n)_{n\in\N}$-driven VARMA model for
$\bepsh$ with real-valued coefficient matrices. It has been shown in
\cite{mokkadem1988mixing}, Theorem 1, that a stable vector ARMA process
is geometrically completely regular provided that the driving noise
sequence is i.i.d. and absolutely continuous with respect to the
Lebesgue measure. A careful analysis of the proof of this result shows
that the existence of an absolutely continuous component of the law of
the driving noise is already sufficient for the conclusion to hold. We
briefly comment on the necessary modifications to the argument. We
first note that under these weaker assumptions, the proof of \cite{mokkadem1988mixing}, Lemma 3, implies that the $n$-step transition
probabilities $P^n(\bx,\cdot)$ of the Markov chain $X$ associated with
a vector ARMA model via its state space representation have an
absolutely continuous component for all $n$ greater than or equal to
some $n_0$. This immediately implies aperiodicity and $\phi
$-irreducibility of $X$, where $\phi$ can be taken as the Lebesgue
measure restricted to the support of the continuous component of
$P^{n_0}(\bx,\cdot)$. The rest of the proof, in particular the
verification of the Foster--Lyapunov drift condition for complete
regularity, is unaltered. This shows that $\bepsh$ is geometrically
completely regular and, in particular, strongly mixing with
exponentially decaying mixing coefficients.
\end{pf*}

\begin{pf*}{Proof of Lemma \ref{MOU}}
By definition, $\M^{(h)}_\nu=\M_\nu(h)$, where $(\M_\nu
(t))_{t\geq0}$
is the solution to
\[
\mathrm{d}\M_\nu(t)=\lambda_\nu\M_\nu(t)\,\mathrm{d}t+\mathrm
{d}\Lb
(t),\qquad\M_\nu(0)=\bzero_m.
\]
Taking the real and imaginary parts of this equation gives
\begin{eqnarray*}
\mathrm{d}\re{\M_\nu(t)} &= & \re{\lambda_\nu\M_\nu(t)}\,
\mathrm
{d}t+\mathrm{d}\Lb(t) = [\re{\lambda_\nu}\re{\M_\nu(t)}-\imag
{\lambda_\nu}\imag{\M_\nu(t)} ]\,\mathrm{d}t+\mathrm{d}\Lb(t),\\
\mathrm{d}\imag{\M_\nu(t)} &= & \imag{\lambda_\nu\M_\nu(t)}\,
\mathrm{d}t
= [\re{\lambda_\nu}\imag{\M_\nu(t)}+\imag{\lambda_\nu}\re{\M
_\nu
(t)} ]\,\mathrm{d}t,
\end{eqnarray*}
and consequently
\[
\mathrm{d} \pmatrix{\displaystyle\re\M_\nu(t)\cr\displaystyle
\imag\M_\nu(t)
}
= [\Lambda_\nu\otimes\I_m ] \pmatrix{\displaystyle\re\M_\nu
(t)\cr\displaystyle\imag\M_\nu(t)
}
\,\mathrm{d}t+ \pmatrix{\displaystyle\I_m \cr\displaystyle 0_{m}
}
\,\mathrm{d}\Lb(t),\qquad\Lambda_\nu= \pmatrix{\displaystyle\re
\lambda_\nu& -\imag\lambda_\nu\cr\displaystyle\imag\lambda_\nu&
\re\lambda_\nu
}
.
\]
Using the fact that $\lambda_\nu\in\R$ for $\nu=1,\ldots,r$ and
$\lambda
_\nu=\overline{\lambda_{\nu+1}}\in\C\backslash\R$ for $\nu
=r+1,r+3,\ldots,N-1$, it follows that $\vecM^{(h)}=\vecM(h)$, where
$(\vecM(t))_{t\geq0}$ satisfies $\mathrm{d}\vecM(t)=G\vecM(t)\,
\mathrm
{d}t+H\,\mathrm{d}\Lb(t)$, and $G=\tilde G\otimes\I_m\in M_{mN}(\R
)$ and
$H=\tilde H\otimes\I_m\in M_{mN,m}$ are given by
\begin{eqnarray*}
\tilde G &=& \operatorname{diag}(\lambda_1,\ldots,\lambda_r,\Lambda
_{r+1},\Lambda_{r+3},\ldots,\Lambda_{N-1}),\\
\tilde H &=& (\underbrace{
\matrix{
1 & \cdots& 1
}
}_{r \mathrm{\ times}}\;\;
\matrix{ 1 & 0 & 1 & 0 & \cdots& 1 & 0
}
)^{\mathrm{T}}.
\end{eqnarray*}
Since $\rank H=m$, the first claim of the lemma is proved. Next, we
show that the controllability matrix $\mathscr{C}\coloneqq\left[
H \enskip GH \enskip\cdots\enskip G^{N-1}H
\right]\in M_{mN}(\R)$ is non-singular. With
$\widetilde{\mathscr
{C}}\coloneqq\left[\vphantom{G^{N-1}\tilde H}
\tilde H \enskip\tilde G\tilde H \enskip\cdots\right.\left. \tilde
G^{N-1}\tilde H
\right]$ and by the properties of the Kronecker product, it follows
that $\mathscr{C}=\widetilde{\mathscr{C}}\otimes\I_m$ and thus
$\det
\mathscr{C} = [\det\widetilde{\mathscr{C}}]^m$. The matrix
$\widetilde
{\mathscr{C}}$ is given explicitly by
\[
\widetilde{\mathscr{C}} =
\pmatrix{
\displaystyle
1 &\lambda_1 & \lambda_1^2 & \cdots & \lambda_1^{N-1}\vspace*{3pt}\cr\displaystyle
\vdots& & & &\vdots\vspace*{3pt}\cr\displaystyle
1 &\lambda_r &\lambda_r^2 &\cdots&\lambda_r^{N-1}\vspace*{3pt}\cr\displaystyle
1 &\re\lambda_{r+1}&\re\lambda_{r+1}^2&\cdots&\re\lambda
_{r+1}^{N-1}\vspace*{3pt}\cr\displaystyle
0 &\imag\lambda_{r+1}&\imag\lambda_{r+1}^2&\cdots&\imag\lambda
_{r+1}^{N-1}\vspace*{3pt}\cr\displaystyle
\vdots& & & &\vdots\vspace*{3pt}\cr\displaystyle
1 &\re\lambda_{N-1}&\re\lambda_{N-1}^2&\cdots&\re\lambda
_{N-1}^{N-1}\vspace*{3pt}\cr\displaystyle
0 &\imag\lambda_{N-1}&\imag\lambda_{N-1}^2&\cdots&\imag\lambda
_{N-1}^{N-1}
}
= T \pmatrix{\displaystyle
1 &\lambda_1 &\lambda_1^2 &\cdots&\lambda_1^{N-1}\vspace*{3pt}\cr\displaystyle
\vdots& & & &\vdots\vspace*{3pt}\cr\displaystyle
1 &\lambda_r &\lambda_r^2 &\cdots&\lambda_r^{N-1}\vspace*{3pt}\cr\displaystyle
1 &\lambda_{r+1}&\lambda_{r+1}^2&\cdots&\lambda_{r+1}^{N-1}\vspace*{3pt}\cr
\displaystyle
\mathrm{i} &\mathrm{i}\overline{\lambda_{r+1}}&\mathrm{i}\overline{\lambda_{r+1}^2}&\cdots
&\mathrm{i}\overline{\lambda_{r+1}^{N-1}}\vspace*{3pt}\cr\displaystyle
\vdots& & & &\vdots\vspace*{3pt}\cr\displaystyle
1 &\lambda_{N-1}&\lambda_{N-1}^2&\cdots&\lambda_{N-1}^{N-1}\vspace*{3pt}\cr
\displaystyle
\mathrm{i} &\mathrm{i}\overline{\lambda_{N-1}}&\mathrm{i}\overline{\lambda_{N-1}^2}&\cdots
&\mathrm{i}\overline{\lambda_{N-1}^{N-1}}
}
\]
with $T\in M_N(\R)$ given by $T=\operatorname{diag} (1,\ldots
,1,R,\ldots,R )$, $R=\frac{1}{2} \left(
{1 \atop-\mathrm{i}}\enskip{ -\mathrm{i} \atop1}
\right)$. Hence, the formula for the determinant of a Vandermonde
matrix (\cite{bernstein2005matrix}, Fact 5.13.3) implies that
\begin{eqnarray*}
&&\det\mathscr{C} = \biggl[(-1)^{({N-r})/{2}}\prod_{1\leq\mu<\nu\leq
r}{(\lambda_\mu-\lambda_\nu)} \mathop{\mathop{\prod}_{\mu,\nu
\in I_{r,N}}}_{ \mu
<\nu}
\imag\lambda_\mu|\lambda_{\mu}-\lambda_{\nu}|^2\\
&&\hphantom{\det\mathscr{C} = \biggl[(-1)^{({N-r})/{2}}\prod_{1\leq\mu<\nu\leq
r}{(\lambda_\mu-\lambda_\nu)} \mathop{\mathop{\prod}_{\mu,\nu
\in I_{r,N}}}_{ \mu
<\nu}}
{}\times|\overline
{\lambda
_\mu}-\lambda_\nu|^2
\mathop{\mathop{\prod}_{1\leq\mu\leq
r}}_{\nu\in
I_{r,N}}{|\lambda_{\mu}-\lambda_{\nu}|^2} \biggr]^m,
\end{eqnarray*}
where $I_{r,N}=\{r+1,r+3,\ldots,N-1\}$. Hence, $\det\mathscr{C}$ is not
zero by Assumption \ref{assumEigen2} and the proof is complete.
\end{pf*}

\section*{Acknowledgements}
This research was supported by the International Graduate School of
Science and Engineering of the Technische Universit\"{a}t M\"{u}nchen
and the TUM Institute for Advanced Study, funded by the German
Excellence Initiative.

\printhistory

\end{document}